\documentclass[oneside, a4paper]{amsart}

	\usepackage[a4paper, margin=1in]{geometry}
	\usepackage{amsmath, amsthm, amssymb, latexsym, mathtools, mathrsfs, eucal}
	\usepackage[all,cmtip]{xy}
	\usepackage{dsfont, soul, stackrel, tikz-cd, graphicx, etoolbox, multirow}
	\usepackage{stmaryrd}
	\DeclareGraphicsExtensions{.pdf,.png,.jpg}
	\usepackage{fancyhdr}
	\usepackage[shortlabels]{enumitem}
	\usepackage[pdfmenubar=true, pdfborder={0 0 0 [3 3]}, colorlinks, citecolor=magenta, urlcolor=blue]{hyperref}
	\numberwithin{equation}{section}
	\setenumerate{listparindent=\parindent}

	\newtheoremstyle{Mytheorem}%
	{1em}{1em}%
	{\slshape}{}%
	{\bfseries}{.}%
	{ }{}

	\newtheoremstyle{Mydefinition}%
	{1em}{1em}%
	{}{}%
	{\bfseries}{.}%
	{ }{}

	\theoremstyle{Mydefinition}
	\newtheorem{statement}{Statement}[section]
	\newtheorem{definition}[statement]{Definition}
	
	\newtheorem{remark}[statement]{Remark}
	
	\newtheorem{example}[statement]{Example}

	\newtheorem*{comment*}{Comment}

	\theoremstyle{Mytheorem}
	\newtheorem{theorem}[statement]{Theorem}
	\newtheorem{corollary}[statement]{Corollary}
	
	\newtheorem{proposition}[statement]{Proposition}
	\newtheorem{lemma}[statement]{Lemma}
	\newtheorem{note}[statement]{Note}

	\newcommand{\G}{GL_2^{+}(\mathbb{Q})}

	\newcommand{\res}{\operatorname{res}}

	\newcommand{\nc}{\newcommand}
	\newcommand{\be}{\begin{eqnarray*}}
	\newcommand{\ee}{\end{eqnarray*}}
	\newcommand{\bea}{\begin{eqnarray}}
	\newcommand{\eea}{\end{eqnarray}}
	\newcommand{\bs}{\begin{split}}
	\newcommand{\es}{\end{split}}
	\newcommand{\bal}{\begin{align}}
	\newcommand{\eal}{\end{align}}

	\nc{\bei}{\begin{itemize}}
	\nc{\eei}{\end{itemize}}
	\nc{\bee}{\begin{enumerate}}
	\nc{\eee}{\end{enumerate}}
	\nc{\bet}{\begin{thm}}
	\nc{\eet}{\end{thm}}
	\nc{\bed}{\begin{defn}}
	\nc{\eed}{\end{defn}}
	\nc{\bel}{\begin{lem}}
	\nc{\eel}{\end{lem}}
	\nc{\bep}{\begin{prop}}
	\nc{\eep}{\end{prop}}
	\nc{\bec}{\begin{corollary}}
	\nc{\eec}{\end{corollary}}
	\nc{\ber}{\begin{rem}}
	\nc{\eer}{\end{rem}}
	\nc{\beex}{\begin{example}}
	\nc{\eeex}{\end{example}}
	\nc{\bpm}{\begin{pmatrix}}
	\nc{\epm}{\end{pmatrix}}
	\nc{\bspm}{\left(\begin{smallmatrix}}
	\nc{\espm}{\end{smallmatrix}\right)}

	\newcommand{\cA}{\mathcal{A}}
	
	\newcommand{\cC}{\mathcal{C}}

	\newcommand{\cO}{\mathcal{O}}

	\newcommand{\cT}{\mathcal{T}}
	\newcommand{\cU}{\mathcal{U}}

	\newcommand{\bC}{\mathbb{C}}

	\newcommand{\BD}{\mathbf{D}}

	\newcommand{\BP}{\mathbf{P}}

	\nc{\frf}{\mathfrak{f}}

	\nc{\frs}{\mathfrak{s}}  
	\nc{\frt}{\mathfrak{t}} 
	\nc{\fru}{\mathfrak{u}}
	\nc{\lsl}{\mathfrak{sl}}
	\nc{\lgl}{\mathfrak{gl}}
	\nc{\upsi}{\underline{\psi}}
	\nc{\uchi}{\underline{\chi}}
	 
	\DeclareMathOperator{\Tr}{Tr}
	\DeclareMathOperator{\PV}{PV}

	\DeclareMathOperator{\Spec}{Spec}

	\DeclareMathOperator{\coker}{coker}

	\DeclareMathOperator{\Cl}{Cl}

	\newcommand{\lra}{\longrightarrow}    
	
	\nc{\surjto}{\twoheadrightarrow}
	\nc{\ts}{\times}
	\nc{\ds}{\displaystyle}
	\nc{\nd}{\noindent}  
	\nc{\ud}{\underline}
	\nc{\ov}{\overline}
	\nc{\maplra}[1]{\buildrel #1 \over \lra}
	\nc{\mapto}[1]{\buildrel #1 \over \to}
	\nc{\setb}[1]{\{  #1\}}
	\nc{\cHom}{\mathcal{H}om}

	\def\a{\alpha}
	\def\b{\beta}

	\def\g{\gamma} \def\G{\Gamma}

	\def\o{\omega} 
	\def\O{\Omega}

	\def\C{\mathbb{C}}
	
	\def\Z{\mathbb{Z}}

\newcommand{\Crit}{\operatorname{Crit}}
\title[Comparison of Frobenius algebra structures on Calabi--Yau toric hypersurfaces]{Comparison of Frobenius algebra structures\\ on Calabi--Yau toric hypersurfaces}

\author{Jeehoon Park}
\address{Jeehoon Park: QSMS, Seoul National University, 1 Gwanak-ro, Gwanak-gu, Seoul 08826, South Korea}
\email{jpark.math@gmail.com}

\author{Philsang Yoo}
\address{Philsang Yoo: Department of Mathematical Sciences \& Research Institute of Mathematics, Seoul National University, 1 Gwanak-ro, Gwanak-gu, Seoul 08826, South Korea}
\email{philsang.yoo@snu.ac.kr}


\begin{document}
\maketitle

\begin{abstract}

%It is well-known that K. Saito's theory of primitive forms and higher residue pairings provides a Frobenius manifold structure on the Jacobian algebra of a homogeneous polynomial with an isolated singularity. This article constructs a Frobenius manifold structure for a class of polynomials that may exhibit non-isolated and non-compact singularities -- specifically, those defining a smooth Calabi--Yau hypersurface in a simplicial Gorenstein toric Fano variety. 

We establish an isomorphism between two Frobenius algebra structures, termed CY and LG, on the primitive cohomology of a smooth Calabi--Yau hypersurface in a simplicial Gorenstein toric Fano variety. As an application of our comparison isomorphism, we observe the existence of a Frobenius manifold structure on a finite-dimensional subalgebra of the Jacobian algebra of a homogeneous polynomial which may exhibit a non-compact singularity locus.

\end{abstract}

\tableofcontents

\section{Introduction}

Let $f   \in \C[\ud z]=\C[z_1,\cdots,z_r] $ be a homogeneous polynomial satisfying $f(0)=0$. Let us assume that $f$ has an isolated singularity at the origin, i.e., the critical locus of $f$
\be
\operatorname{Crit}(f):=\left \lbrace  (a_1,\cdots,a_r) \in \C^r \;\middle|\; \frac{\partial  f}{\partial z_i}(a_1,\cdots,a_r) = 0 \text{ for }  i =1, \cdots, r\right\rbrace 
\ee
is a singleton consisting of the origin. Then the Jacobian algebra,
\be
R(f):=\C[\ud z]/ J(f) \quad \text{where}\quad   J(f)=\left\langle  \frac{\partial f}{\partial z_1} , \cdots , \frac{\partial f}{\partial z_r} \right\rangle,
\ee
is a finite-dimensional vector space over $\C$, and it is well-known that K. Saito's theory of primitive forms and higher residue pairings (\cite{Saito}, \cite{Saitoh}, and \cite{ST}) provides a (formal) non-trivial Frobenius manifold structure on the Jacobian algebra $R(f)$.\footnote{In fact, Saito's construction yields a stronger result: an analytic Frobenius manifold structure, rather than just a formal one. This is achieved by solving a Birkhoff factorization problem.}

A natural question is whether a Frobenius manifold can be associated with a polynomial $f$ that has a non-isolated singularity. Clearly, Saito's theory cannot be applied; after all, $R(f)$ is not even finite-dimensional unless the critical locus consists of finitely many points. Earlier works \cite{Sa}, \cite{OV} explored the case of non-isolated critical loci using characteristic $p$ methods. More recently, Li and Wen \cite{LW} introduced Hodge-theoretic methods to construct a Frobenius manifold for polynomials with compact critical loci, drawing parallels to the work of Barannikov and Kontsevich \cite{BK}. These approaches overlap in some cases but also address scenarios not covered by others.
% Recently, Li and Wen \cite{LW} introduced Hodge-theoretic methods to construct a Frobenius manifold for polynomials with compact critical loci, drawing parallels to those of Barannikov and Kontsevich \cite{BK}. There are also works on non-isolated critical loci case using characteristic $p$ method \cite{Sa}, \cite{OV}; these works have cases in common and also cases not covered by others. 
In this article, through a comparison isomorphism between two Frobenius algebra structures (Theorem \ref{mto}) we provide an indirect way of showing the existence of Frobenius manifold structures for a class of polynomials with possibly non-compact critical loci.

To proceed and find a promising direction to pursue, let us revisit the case of an isolated singularity. Suppose $f$ is a homogeneous polynomial of degree $r$ that defines a smooth Calabi--Yau hypersurface $X_f$ embedded in the projective space $\BP^{r-1}$, yielding a compact K\"ahler manifold $X:=X_f(\C)$ of dimension $r-2$. According to Griffith's theorem \cite{Gr69}, there is a $\C$-vector space isomorphism\footnote{Without Calabi--Yau condition, one still has the following isomorphism
\be
\phi: \bigoplus_{a=0}^{\infty} R(f)_{(a+1)\deg f - r} \cong H_{\text{pr}}^{r-2}(X).
\ee}
\bea \label{Gr}
\phi: \bigoplus_{a=0}^{\infty} R(f)_{a r} \cong  H_{\text{pr}}^{r-2}(X)
\eea
where $R(f)_{k}=\dfrac{\C[\ud z]_{k}}{J(f) \cap \C[\ud z]_{k}}$ is the homogeneous component of $R(f)$ of degree $k$ and $H_{\text{pr}}^{r-2}(X) = H_{\text{pr}}^{r-2}(X, \C)$ is the primitive middle-dimensional cohomology of $X$ with the coefficient $\C$. Moreover, if we consider the Hodge decomposition $H_{\text{pr}}^{r-2}(X) \cong  \bigoplus_{a=0}^{r-2} H_{\text{pr}}^{r-2-a,a}(X)$, then $\phi$ maps $R(f)_{ar}$ to $H_{\text{pr}}^{r-2-a,a}(X)$ for each $a$. This is compatible with the Macaulay theorem \cite{Mac} that $R(f)_{ar}=0$ for $a \geq r-1$.
%In fact, it is easy to check that Saito's construction restricts to $A(f)$ and provides a (formal) non-trivial Frobenius manifold structure on $A(f)$.

In fact, we can choose a $\C$-vector space isomorphism $\phi$ so that it preserves additional structures. Specifically, there exists a trace map $\Tr\colon R(f)_{(r-2)r}\cong  \C$ corresponding to the integration map $ \int_X \colon  H_{\text{pr}}^{r-2,r-2}(X)\cong  \C $ in that  the following diagram
\[\begin{tikzcd}[row sep=0.5em, column sep=5em]
{R(f)_{ar} \times R(f)_{br}} & {R(f)_{(r-2)r}} \\ && \C \\ {H_{\text{pr}}^{r-2-a,a}(X) \times H_{\text{pr}}^{r-2-b,b}(X)} & {H_{\text{pr}}^{r-2,r-2}(X)} \arrow["{\text{mul}}", from=1-1, to=1-2] \arrow["{\phi\times \phi}", from=1-1, to=3-1] \arrow["\Tr", from=1-2, to=2-3]  \arrow["\wedge", from=3-1, to=3-2] \arrow["{\int_X}"', from=3-2, to=2-3]
\end{tikzcd}\]
commutes \cite[Theorem 3]{CG}, where 
\be
\text{mul}(u,v) = (-1)^b u \cdot v \qquad\text{for}\quad  u \in R(f)_{ar},\ v \in R(f)_{br}.
\ee
% up to sign, i.e.,
%\be
%Tr(u \cdot v)= (-1)^{b} \int_X \phi(u) \wedge \phi(v), \quad u \in R(f)_{ar}, v \in R(f)_{br}
%\ee
whenever $a+b =r-2$ (see the proof of Proposition \ref{mpro} given at the end of Subsection \ref{ss3.2} for details for where the sign $(-1)^b$ comes from).  Note that the bilinear map $\text{mul}$ is symmetric (respectively, skew-symmetric) if $r-2$ is even (respectively, $r-2$ is odd).

Let us denote the sub-algebra $\bigoplus_{a=0}^{\infty} R(f)_{a r}$ of $R(f)$ by \[A(f):=\bigoplus_{a=0}^{\infty} R(f)_{a r}\] and the above symmetric bilinear pairing (by declaring $\langle u, v\rangle =0$ for $u\in R(f)_{ar}$, $v \in R(f)_{br}$ unless $a+b=r-2$) by
\be
\langle - , - \rangle: A(f) \times A(f) \to \C.
\ee
The bilinear pairing extends to a symmetric bilinear pairing $\langle -, - \rangle:R(f) \times R(f) \to \C$ given by $\langle u,v \rangle =  \Tr(u\cdot  v)$; we still have $\langle u, v\rangle = 0$ unless $\deg u+ \deg v= (r-2)r$. Then the triple $(R(f), \cdot, \langle - , - \rangle)$ is a Frobenius algebra and $(A(f), \cdot, \langle - , - \rangle)$ is its Frobenius subalgebra. As the aforementioned theory of K. Saito promotes the Frobenius algebra $(R(f),\cdot, \langle - , - \rangle)$ to a Frobenius manifold structure on $R(f)$, it restricts to $(A(f), \cdot, \langle - , - \rangle)$ and provides a formal Frobenius manifold structure on $A(f)$ as well.

%in the following sense. Let $\ud t=\{t^\a \}$ be a coordinate (around 0) of the universal unfolding parameter space $R(f)\cong \C^\nu$ which defines a $\C$-basis $\{ u_\a \} $ of $R(f)$. Then Saito's theorem says that 
%\begin{theorem}\label{Saito} (\cite{Saito}, \cite{Saitoh})
%There exist a 3-tensor $A_{\a\b}^\g (\ud t) \in \C \llbracket  \ud t \rrbracket  $ and a 2-tensor $g_{\a\b} \in \C$, which satisfy the axioms for Frobenius manifold structure (see Definition \ref{Frob}), where
%\be
%u_\a \cdot u_\b = \sum_{\g} A_{\a\b}^\g (\ud t) u_\g
%\ee
%defines a structure of commutative $\C\llbracket \ud t \rrbracket$-algebra on $R(f) \otimes_\C \C\llbracket \ud t \rrbracket$, and
%\be
%\langle u_\a, u_\b \rangle = g_{\a\b}
%\ee 
%defines a symmetric $\C\llbracket \ud t \rrbracket$-linear pairing on $R(f) \otimes_\C \C\llbracket \ud t \rrbracket$. Moreover, $A_{\a\b}^{\g}(0)$ and $g_{\a\b}$ gives the initial Frobenius algebra $(R(f), \langle - , - \rangle, \Tr)$. 
%\end{theorem}

The goal of this article is to identify a class of polynomial functions $f$ on $\C^r$ for which the critical locus $\Crit(f)$ may be non-isolated and non-compact in $\C^r$, yet the following are satisfied:
\begin{enumerate} 
\item there exists a finite-dimensional $\C$-vector subspace $A(f) \subset R(f)$ that is equipped with a Frobenius algebra structure;
\item there exists a formal non-trivial Frobenius manifold structure on $A(f)$ that extends the Frobenius algebra $A(f)$.
\end{enumerate}

%there still exists a formal non-trivial Frobenius manifold structure on a certain finite-dimensional (over $\C$) subalgebra of the Jacobian algebra $R(f)$ (note that $R(f)$ is possibly infinite-dimensional over $\C$).

In the rest of the introduction, we sketch our approach. Our setup begins with an $m$-dimensional simplicial Gorenstein toric Fano variety $\BP$ whose toric homogeneous coordinate ring is $\C[\ud z] =\C[z_1,\cdots,z_r] $. Thus $\BP$ is a complete normal orbifold whose anti-canonical divisor is Cartier and ample.\footnote{These assumptions on $\BP$ are crucial in using the Batyrev--Cox theorem \cite[Theorem 10.13]{BatCox}; see Theorem \ref{thm:BatCox}.} Suppose $f \in \C[\ud z]$ defines a smooth ample Calabi--Yau hypersurface $X_f$ in $\BP$ and a compact K\"ahler manifold $X:=X_f(\C)$. In this case, the coordinate ring $\C[\ud z]$ is graded by the class group of $\BP$. As we write $R(f)$ for the Jacobian algebra, for an element $\a$ in the class group of $\BP$, we denote the homogeneous component of $R(f)$ with degree $\a$ by $R(f)_{\a}$. If we define \[A(f)=\bigoplus_{a=0}^{\infty} R(f)_{a \b},\] where $\b$ is the anti-canonical class of $\BP$, then $A(f)$ is still a subalgebra of $R(f)$, whose multiplication we denote by $\bullet_{\text{LG}}$. 

Now our main question is whether one can promote this algebra structure on $A(f)$ to a Frobenius algebra structure and moreover construct a non-trivial Frobenius manifold structure. The main difficulty (in the non-isolated and non-compact case) in constructing a Frobenius manifold structure on $A(f)$ lies in the lack of the higher residue pairing associated to $A(f)$ and its primitive forms; one can not apply \cite{LLS}. Since  $\Crit(f)$ is non-compact, it is not possible to apply the $L^2$-Hodge theoretic method of \cite{LW}, either.

Our approach is 
\begin{enumerate}
	\item to recall Batyrev--Cox's theorem \cite[Theorem 10.13]{BatCox} 
\be
A(f) \cong H_{\text{pr}}^{m-1}(X),
\ee
which generalizes \eqref{Gr}\footnote{When $\BP=\BP^{r-1}$, the anti-canonical class $\beta$ of $\BP^{r-1}$ is $\mathcal O_{\BP^{r-1} }(r)$ which corresponds to the integer $r$ under the isomorphism $\Cl (\BP^{r-1})\cong \Z $.} (see Theorem \ref{thm:BatCox}) and induce a trace pairing on $A(f)$ via this isomorphism (see Proposition \ref{mpro}); this yields what we call the LG (Landau--Ginzburg) Frobenius algebra and denote by $(A(f),\bullet_{\text{LG}} , \langle -, - \rangle_{\text{LG}} ) $, and equivalently, $(H^{m-1}_{\text{pr}}(X),\bullet_{\text{LG}} , \langle -, - \rangle_{\text{LG}} ) $;
	\item to identify it with the CY (Calabi--Yau) Frobenius algebra structure \[  (H_{\text{pr}}^{m-1}(X), \bullet_{\text{CY}},  \langle -,-\rangle_{\text{CY} }  )\] on $H_{\text{pr}}^{m-1}(X)$ (see Definition \ref{CYFrob}), which has been proved to extend to a formal Frobenius manifold structure on $H_{\text{pr}}^{m-1}(X)$ by Barannikov and Kontsevich \cite{BK}.
\end{enumerate}

More precisely, we show that two (CY and LG) Frobenius algebra structures are equivalent and simply transport the Frobenius structure on the CY Frobenius algebra to the LG Frobenius algebra.

\begin{theorem}\label{mto}
Let $\BP$ be an \( m \)-dimensional simplicial Gorenstein toric Fano variety with a toric homogeneous coordinate ring \( S = \C[z_1, \ldots, z_r] \). Let \( f \in S \) be a polynomial of degree given by the anti-canonical divisor \( \b \) of \( \BP \), defining a smooth Calabi--Yau hypersurface \( X_f \) in \( \BP \). Assume that the map 
\bea 
\label{extraisom} H^{m-2}(\BP) \xrightarrow{\cup [X]} H^{m}(\BP) 
\eea
 is an isomorphism, where \( [X] \in H^2(\BP) \) denotes the cohomology class of \( X = X_f(\C) \). Then the Landau--Ginzburg Frobenius algebra \( ( H^{m-1}_{\text{pr}}(X) ,\bullet_{\text{LG}} , \langle -, - \rangle_{\text{LG}} ) \) is isomorphic to the Calabi--Yau Frobenius algebra \( (H_{\text{pr}}^{m-1}(X), \bullet_{\text{CY}}, \langle -,-\rangle_{\text{CY}} ) \).
\end{theorem}

Based on Theorem \ref{mto}, we can transport the formal Frobenius manifold structure on the CY side to the LG side and deduce the existence of a formal Frobenius manifold structure on $A(f)$ extending $(A(f),\bullet_{\text{LG}} , \langle -, - \rangle_{\text{LG}} ) \cong (H_{\text{pr}}^{m-1}(X), \bullet_{\text{CY}},  \langle -,-\rangle_{\text{CY} }  )$:

\begin{corollary} \label{mtt}
There exists a formal non-trivial Frobenius manifold structure on $A(f) $ which extends the Frobenius algebra $(A(f),\bullet_{\text{LG}} , \langle -, - \rangle_{\text{LG}} )$.
\end{corollary}

We briefly explain the structure of the article. In Section \ref{sec2}, we explain two Frobenius algebra structures on $H_{\text{pr}}^{m-1}(X)$: the CY Frobenius algebra structure in Subsection \ref{ss2.1} and the LG Frobenius algebra structure in Subsection \ref{ss2.2}. In Section \ref{sec3}, we compare these two Frobenius algebra structures (Theorem \ref{mto}) and provide its application (Corollary \ref{mtt}). In Subsection \ref{ss3.1}, we present the crucial result required for the proof of Theorem \ref{mto}. In Subsection \ref{ss3.2}, we prove main results. In Subsection \ref{ss3.3}, we deduce Corollary \ref{mtt}. Finally, in Subsection \ref{ss3.4}, we compare it with the case of isolated singularities and provide a nontrivial example in which $\BP$ has Picard rank 2 and $\Crit(f)$ is non-isolated and non-compact.
%
%how we enhance the Batyrev--Cox isomorphism to a ring isomorphism.
%Subsection \ref{ss2.1} is devoted to explain the CY ring structure on $H_{\text{pr}}^{m-1}(X)$.
%In Subsection \ref{ss2.2}, we explain the LG ring structure on $H_{\text{pr}}^{m-1}(X)$ and state the main result (Theorem \ref{mthm}). We prove Theorem \ref{mthm} by applying Loyola's toric Carlson--Griffiths theorem in Subsection \ref{ss2.3}.
%Section \ref{sec3} is devoted to prove Theorem \ref{mtt} by applying the results in Section \ref{sec2}. In Subsection \ref{ss3.1}, we review the CY Frobenius manifold construction of Barannikov--Kontsevich. In Subsection \ref{ss3.2}, we construct the LG Frobenius algebra and prove Theorem \ref{mtt}.

\vspace{1em}

\textbf{Acknowledgement}:
J.P. would like to thank R. Villaflor Loyola for a detailed explanation for Lemma \ref{mac}.
J.P. also thanks V. Batyrev for answering a question on the toric Macaulay theorem. P.Y. thanks S. Li for valuable discussions on a related topic.

%The work of Jeehoon Park was partially supported by BRL (Basic Research Lab) through the NRF (National Research Foundation) of South Korea (NRF-2018R1A4A1023590). 
%The work of Jeehoon Park was supported by the National Research Foundation (NRF) grant of South Korea (NRF- 2018R1A4A1023590, No.2020R1A5A1016126, and NRF-2021R1A2C1006696) funded by the Korea government (MSIT). 
%The author was supported by Samsung Science \& Technology Foundation (SSTF-BA1502),
%the National Research Foundation of Korea (NRF-2021R1A2C1006696) and the National Research Foundation of Korea (NRF) grant funded by the Korea government (MSIT) (No.2020R1A5A1016126).
Jeehoon Park was supported by the National Research Foundation of Korea (NRF-2021R1A2C1006696) and the National Research Foundation of Korea (NRF) grant funded by the Korea government (MSIT) (No.2020R1A5A1016126). Philsang Yoo was supported by the National Research Foundation of Korea (NRF) grant funded by the Korea government(MSIT) (No. 2022R1F1A107114212).

\section{Two Frobenius algebra structures}\label{sec2}

In this section, we will find two Frobenius algebra structures (namely, CY Frobenius algebra and LG Frobenius algebra) on the primitive middle-dimensional cohomology $H_{\text{pr}}^{m-1}(X)$. This will be done by transport of structures through the following isomorphisms of graded vector spaces
\[
\begin{tikzcd}
H_{\text{pr}}^0 (\PV(X)) & H_{\text{pr}}^{m-1}(X) \arrow[l, "\cong" above] & A(f) .\arrow[l, "\cong" above]
\end{tikzcd}
\]

The CY Frobenius algebra structure exists on the cohomology of the space $\PV(X)$ of polyvector fields for any compact Calabi--Yau manifold $X$ by the work of Barannikov--Kontsevich in \cite{BK}, which we briefly recall in Subsection \ref{ss2.1}. Then, in Subsection \ref{ss2.2}, we will construct the LG Frobenius algebra structure on $H_{\text{pr}}^{m-1}(X)$; this is based on the Batyrev--Cox theorem, which says that $H_{\text{pr}}^{m-1}(X)$ is isomorphic to a certain subring $A(f)$ of the Jacobian ring $R(f)=S/J(f)$ where $J(f)$ is the Jacobian ideal of $f$.

%By contracting with a nowhere-vanishing holomorphic form, which exists due to the Calabi--Yau assumption of $X_f$, one can equip $H_{\text{pr}}^{m-1}(X)$ with a ring structure (called the Calabi--Yau ring structure) coming from the product of polyvector fields on $X_f$. 

%a global section of $\Omega^{n+1}_{\BP^{n+1}}(n+2) \simeq \cO_{\BP^{n+1}}$:
%\begin{align*}
%\Omega = \sum_{i=0}^{n+1} (-1)^i x_i (dx_0\wedge \cdots \wedge \hat {d x_i} \wedge \cdots \wedge dx_{n+1}),
%\end{align*}
%which is unique up to constant.

%\subsection{dGBV algebra of Barannikov--Kontsevich}
%\subsection{CY ring structure}\label{ss2.1}	

\subsection{CY Frobenius algebra}\label{ss2.1}

Let $X$ be a compact Calabi--Yau manifold of dimension $n$. We fix a choice of a nowhere-vanishing holomorphic $n$-form $\Omega_X$ on $X$. %Using this $\Omega_X$, there is a way (due to \cite{BK}) to put a ring structure on the middle-dimensional singular cohomology $H^n(X, \C)$; this ring structure plays an important role on the $B$-side of the mirror symmetry conjecture as is used in \cite{BK} by Barannikov and Kontsevich.
Let $\cT_X$ (respectively, $\cT^*_X$) be the holomorphic tangent (respectively, cotangent) bundle on $X$. We write $\O_X^k:=\wedge^k \cT^*_X$.

Let $\PV^{i,j}(X):=\cA^{0,j}(X, \wedge^i \cT_X)$ denote the space of smooth $(0,j)$-forms valued in $\wedge^i \cT_X$. Consider the isomorphism $\cA^{0,j}(X, \wedge^i \cT_X) \cong \cA^{0,j}(X, \Omega^{n-i}_X) $ given by contracting with $\Omega_X$, which we denote by $\vdash \Omega_X $. This induces an isomorphism
\be
 c_{\O_X}: \bigoplus_{i,j=0}^n \cA^{0,j}(X, \wedge^i \cT_X) \xrightarrow{\cong} \bigoplus_{i,j=0}^n\cA^{0,j}(X, \Omega^{n-i}_X)
\ee
and hence an isomorphism $H^j(X , \wedge^i \cT_X )\cong H^{n-i,j}(X)$ for each $0\leq i,j\leq n$ after taking $\overline{\partial }$-cohomology. Moreover, by the Hodge decomposition $H^k(X) \cong \bigoplus_{p+q=k} H^{p,q}(X)$ and by the degeneration of the Hodge-to-de Rham spectral sequence, this induces an isomorphism $H^k (X) \cong  \bigoplus_{j-i=k-n} H^j(X , \wedge^i \cT_X ) $ for each $0\leq k\leq n$. %Henceforth we write $H^{a,b}(\PV(X))= H^b(X,\wedge^a \cT_X)$.

Note that the space $ \bigoplus_{i,j=0}^n \PV^{i,j}(X)$ has a product structure coming from the product of holomorphic polyvector fields valued in anti-holomorphic differential forms:
\be \label{CYR}
\wedge \colon  \cA^{0,j}(X, \wedge^i \cT_X)\times \cA^{0,l}(X, \wedge^k \cT_X) \to \cA^{0,j+l}(X, \wedge^{i+k} \cT_X).
\ee
Now let us construct a Frobenius algebra that will correspond to $H^n(X)$. Let us set
\be \label{pvgrading}
\PV^r(X):= \bigoplus_{r=j-i}\PV^{i,j}(X).
\ee
Then in the case of $r=0$, the contraction map $ c_{\O_X}$ restricts to an isomorphism
\be
 c_{\O_X}: \bigoplus_{j=0}^n \cA^{0,j}(X, \wedge^{j} \cT_X) \xrightarrow{\cong} \bigoplus_{j=0}^n\cA^{0,j}(X, \Omega^{n-j}_X)
\ee
and induces an isomorphism $ H^0(\PV(X) ):=\bigoplus_{j=0}^n  H^{j,j}(\PV(X))\xrightarrow{\cong} H^n(X)$ where $ H^{j,j}(\PV(X)):=H^j(X, \wedge^j \cT_X)$. Under this isomorphism $H^{j,j}(\PV(X))$ corresponds to $H^{n-j,j}(X)$.
%$H^{0,0}(\PV(X)) \simeq \bigoplus_{j=0}^n H^j(X, \wedge^j \cT_X) \xrightarrow{\cong} H^n(X)$
Also, it has a trace map $\Tr_{\text{CY}}: \PV(X) \to  \C$ given by
\bea\label{tracecy}
\Tr_{\text{CY}}(\o) = \int_X \left(\o \vdash \O_X\right) \wedge \O_X.
\eea
Clearly, the trace map vanishes unless $\omega \in \PV^{n,n}(X)$ and in fact it induces an isomorphism 
$ H^{n,n}(\PV(X))  \cong \C$. 
%$ H^n(X, \wedge^n \cT_X):=\PV^{n,n}(X)/\bar\partial(\PV^{n,n-1}(X))  \cong \C$.
Using $\Tr_{\text{CY}}$, we define a symmetric bilinear pairing $\langle -,- \rangle_{\text{CY}}$ on $H^0(\PV(X))$ by
\bea \label{cypair}
\langle \o, \eta \rangle_{\text{CY}} :=\Tr_{\text{CY}}(w\wedge  \eta)= \int_X \left((\o \wedge  \eta ) \vdash \O_X\right) \wedge \O_X,
\eea
for $\o, \eta \in \PV(X)$. This yields a Frobenius algebra $( \bigoplus_{i,j=0}^n H^j(X , \wedge^i \cT_X )  , \wedge ,  \langle  -,- \rangle_{\text{CY}})$.

By construction, the product map of the form
\[\wedge \colon  \cA^{0,i}(X, \wedge^i \cT_X)\times \cA^{0,j}(X, \wedge^j \cT_X) \to \cA^{0,i+j}(X, \wedge^{i+j} \cT_X)
\] is closed in $H^0(\PV(X) )$ and the trace map $\Tr_{\text{CY}} \colon \PV(X)\to \C $ is still nontrivial on $H^0(\PV(X))$ as the canonical map $ H^{n,n}(\PV(X))\to H^0(\PV(X))$ is injective. Hence we have a nontrivial Frobenius algebra $(H^0 (\PV(X)), \wedge ,  \langle  -,- \rangle_{\text{CY}})$.

\begin{definition}\label{CYFrob}
 We call the triple $(H ^0(\PV(X)), \wedge,  \langle  -,- \rangle_{\text{CY}})$ the \emph{CY Frobenius algebra} of a compact Calabi--Yau manifold $X$ of dimension $n$. By transport of structure via $c_{\Omega_X}$, it defines the \emph{CY Frobenius algebra} $(H^n(X), \bullet_{\text{CY}},  \langle -,-\rangle_{\text{CY} } )$.
\end{definition}

We make the following observation.

\begin{proposition} \label{prop:CY trace}
The trace map $\Tr_{\text{CY}}$ in \eqref{tracecy} makes the following diagram
\[\begin{tikzcd}[row sep=0.5em, column sep=5em]
	H^{0}(\PV(X)) \times H^{0}(\PV(X)) & H^{0}(\PV(X)) \\
	&& \mathbb{C} \\
	{H^{n}(X) \times H^{n}(X)} & {H^{2n}(X)}
	\arrow["{\wedge^*}", from=1-1, to=1-2]
	\arrow["{c_{\O_X} \times c_{\O_X}}", from=1-1, to=3-1, swap]
	\arrow["\Tr_{\text{CY}}", from=1-2, to=2-3]
	\arrow["\cup", from=3-1, to=3-2]
	\arrow["{\int_X}"', from=3-2, to=2-3]
\end{tikzcd}
\]
 commutes, where
 \be
  \omega \wedge^* \eta := (-1)^b \omega \wedge \eta \qquad \text{for}\quad  \o \in  \PV^{a,a}(X),\ \eta \in \PV^{b,b}(X)
 \ee
 whenever $a+b =n$.\footnote{Note that the bilinear map $\wedge^*$ is symmetric (respectively, skew-symmetric) if $n$ is even (respectively, $n$ is odd).}
 % 0 \leq a, b \leq n$.
%  More precisely we have 
%\be 
%\int_X c_{\Omega}(\omega \wedge \eta) \wedge \O_X = (-1)^{ b } \int_X c_{\Omega}(\omega) \cup c_{\Omega}(\eta)
%\ee
%or $\o \in  \PV^{n-b,n-b}(X)$ and $\eta \in \PV^{b,b}(X)$.
 %\footnote{Here the wedge product is commutative, whereas the cup product is graded-commutative. In particular, when $n$ is even, the diagram commutes with no ambiguity. When $n$ is odd, the sign ambiguity is unavoidable.  }
\end{proposition}

 \begin{proof}
One can check that
\be 
\int_X \left( (\o \wedge  \eta) \vdash \O_X\right) \wedge \O_X = (-1)^{ nb + (n+1)a } \int_X (\o \vdash \O_X) \cup (\eta \vdash \O_X)
\ee
for $\o \in  \PV^{n-a,n-b}(X)$ and $\eta \in \PV^{a,b}(X)$.
If we put $a=b$, then we get the desired result:
\be
\Tr_{\text{CY}}(\o \wedge^* \eta) = \int_X c_{\Omega_X}(\o) \cup c_{\Omega_X}(\eta)\qquad \text{for}\quad  \o \in  \PV^{n-b,n-b}(X),\ \eta \in \PV^{b,b}(X).
\ee
 \end{proof}

\subsection{LG Frobenius algebra}\label{ss2.2}

In this subsection, we discuss a Frobenius algebra structure on the middle-dimensional primitive cohomology of a class of Calabi--Yau manifolds. This structure is induced from a subalgebra of the Jacobian algebra and is therefore referred to as the LG Frobenius algebra structure. Much of the content here reorganizes materials from \cite{BatCox}, which should be consulted for additional details and broader context.

Let  $\BP$ be an $m$-dimensional simplicial Gorenstein toric Fano variety, which by definition is a projective orbifold. More concretely, there exists a simplicial fan $\Sigma$ so that $\BP=\BP_\Sigma$ is the associated toric variety. Moreover, we have the followings:
\begin{itemize}
	\item There exists a homogeneous coordinate ring $S=S(\Sigma)=\C[z_1, \cdots, z_r]$ of $\BP_\Sigma $ graded by the class group $\Cl(\BP_\Sigma )$ whose rank is $r-m$. Here $r:={\#} \Sigma(1)$ is the number of one-dimensional rays of $\Sigma$. Each variable $z_i$ corresponds to the torus-invariant divisor $D_i$, which is associated to a ray in $\Sigma(1)$ generated by a primitive element $\rho_i$.
	\item The anti-canonical divisor $\b=-K_{\BP_\Sigma}=\sum_{\rho \in \Sigma(1)} D_\rho \in \Cl(\BP_\Sigma)$ is Cartier and ample.
\end{itemize}

Let $X_f \subset \BP$ be a quasi-smooth Calabi--Yau hypersurface (of dimension $m-1$) defined by $f \in S $ of degree $\deg f =\b=\sum_{\rho \in \Sigma(1)} D_\rho$ so that $X_f$ is ample.
 %$\BD(\Simga):=\Spec \C[\Cl(BP)]$
%Let
%\be
%A_\b:=\bigoplus_{a=0}^{\infty}S_{a\b} \subset S
%\ee
%denote the subring of $S$.
Recall that $J(f)$ is the Jacobian ideal of $S$ generated by $f_1, \cdots, f_{r}$, where $f_i =\frac{\partial f}{\partial z_i}$, and we set
\be
R(f):=S/J(f) \qquad \text{and}\qquad  A(f):= \bigoplus_{a=0}^\infty R(f)_{a\b}.
\ee
Then $A(f)$ is a subalgebra of $R(f)$. This gives an algebra structure on $A(f)$ that we denote by $\bullet_{\text{LG}}$. On the other hand, there is no canonical trace map on $A(f)$ because there is none on the Jacobian algebra $R(f)$ in this generality. In order  to find  $\Tr_{\text{LG}} \colon A(f)\to \C $ in an analogous way to Proposition \ref{prop:CY trace}, we first turn to finding a graded vector space that is isomorphic to $A(f)$.

\begin{definition}
The \emph{primitive middle-dimensional cohomology} of $X=X_f(\C)$ is defined by
 \be
 H^{m-1}_{\text{pr}}(X) := \coker \left(H^{m-1}(\BP, \C) \xrightarrow{i^*}  H^{m-1}(X, \bC) \right)
 \ee
 where $i: X \to  \BP$ is the canonical embedding.
 \end{definition}
 
 \begin{note}
 The CY Frobenius structures on $H^0(\PV(X))\cong H^n(X)$ in Definition \ref{CYFrob} induces CY Frobenius algebra structures on $H^0_{\text{pr}}(\PV(X)) \cong  H^{n}_{\text{pr}}(X)$. Henceforth, the CY Frobenius structure refers to the induced one on the primitive cohomology.
 \end{note}

Recall that for a normal variety $Y$, the sheaf $\Omega_Y^{k}$ of $k$-forms is not necessarily well-behaved. On the other hand, the sheaf $\hat \O_Y^k := (\Omega_Y^{k})^{\vee \vee}$ is better behaved with which much of general theory is concretely developed in a toric setting. We refer to \cite{BatCox, CLS} for further discussion. We fix a generator of $H^0(\BP_\Sigma,  \hat \Omega_{\BP_\Sigma}^m(\b))$ 
where $\b$ is the anti-canonical class of $\BP_\Sigma$:
\bea \label{vf}
\O := \sum_{|I|=m} \det(\rho_I) \widehat {z_{I}} d z_I
\eea
where for  $I=(i_1, \cdots, i_m)$ we write $d z_I = dz_{i_1} \wedge \cdots \wedge dz_{i_m}$, $\widehat {z_{I}} =\prod_{i \notin I}z_i$, and $\det(\rho_{I})=\det (\rho_{i_1} | \cdots | \rho_{i_m})$ with $\rho_i$ primitive generators of the rays of $\Sigma$. Batyrev and Cox showed the following theorem. 
 \begin{theorem}[Batyrev--Cox] \cite[Theorem 10.13]{BatCox}  \label{thm:BatCox}
 Let $X$ be a quasi-smooth ample Calabi--Yau hypersurface in an $m$-dimensional simplicial Gorenstein toric Fano variety $\BP$. Assume \eqref{extraisom}.\footnote{If $m$ is even, then, for $a = \frac{m}{2}$, the map $\phi_\Omega \colon R(f)_{a\b} \to H^a_{\text{pr}}(X, \hat  \O_X^{m-1-a})$ might not be an isomorphism without the assumption \eqref{extraisom}.}
 Then there is a canonical isomorphism of vector spaces
 \[ \phi_\Omega \colon  R(f)_{a \b}\to H^a_{\text{pr}}(X, \hat \O_X^{m-1-a})\]
for each $0\leq a\leq m-1$, defined by
 \be
 \phi_\O ([F(\ud z)]) = \res \left((-1)^a a!  \frac{F \Omega}{ f^{a+1}} \right) \in H^a_{\text{pr}}(X,  \hat \O_X^{m-1-a}) \qquad \text{for}\qquad F(\ud z) \in R(f)_{a\b}.
 \ee
Moreover, $H^{m-1}_{\text{pr}}(X, \C)$ has a pure Hodge structure and hence there is a decomposition $H^{m-1}_{\text{pr}}(X, \C)\cong \bigoplus_{a=0}^{m-1}H^a_{\text{pr}}(X, \hat \O_X^{m-1-a})$. Therefore, there is a graded vector space isomorphism (depending on $\Omega$)
 \be \label{BCe}
 \phi_\O:  \bigoplus_{a=0}^{m-1} R(f)_{a \b}   \xrightarrow{\cong} H^{m-1}_{\text{pr}}(X, \C).
 \ee
 \end{theorem}
 %Since $A/Jac(f) \cap A$ has the ring structure induced from $\C[\ud x]/Jac(f)$, the isomorphism $\phi_\O$ puts the ring structure on $H^n_{\text{pr}}(X, \C)$, which we call \textit{the Griffiths ring structure} or \textit{the LG (Landau--Ginzburg) ring structure} on $H^n_{\text{pr}}(X,\C)$.

Given that $H^{m-1}_{\text{pr}}(X, \C)\cong \bigoplus_{a=0}^{\infty}H^a_{\text{pr}}(X, \hat \O_X^{m-1-a})$ holds for a degree reason, it is natural to ask whether $ A(f)=\bigoplus_{a=0}^{\infty }R(f)_{a \b}$ also coincides with $\bigoplus_{a=0}^{m-1} R(f)_{a \b}$. It follows from proving the following lemma.

\begin{lemma}\label{mac}
Let $X$ be a quasi-smooth Calabi--Yau hypersurface in an $m$-dimensional simplicial Gorenstein toric Fano variety $\BP$, defined by $f\in S$ of degree $\deg f = \b =-K_{\BP}=\sum_{\rho \in \Sigma(1)} D_\rho$. Then we have
\be \label{must}
R(f)_{p \beta} = 0 \qquad\text{for}\qquad  p \geq m.
\ee
\end{lemma}
\begin{proof}
By \cite[Corollary 10.2]{BatCox}, we have a $\C$-vector space isomorphism
\bea \label{bca}
H^m(\BP \setminus X)\cong   \bigoplus_{p=0}^{m} \frac{H^0(\BP,\hat  \Omega_\BP^m ((p+1)X)) }{H^0(\BP,\hat  \Omega_\BP^m (pX)) + dH^0(\BP,\hat  \Omega_\BP^{m-1} (pX))}.
\eea
By \cite[Theorem 9.7]{BatCox},  we have
\be
H^0(\BP,\hat  \O_{\BP}^m ((p+1) X)) =\left\{ \frac{ u \O}{f^{p+1}} :  u \in S_{p \b}\right\}
\ee
for any $p \geq 0$, and the map 
\be
\psi: H^0(\BP, \hat  \O_{\BP}^m ((p+1) X)) \xrightarrow{\cong}  S_{p\b}, \quad p \geq 0
\ee
defined by $\psi (\frac{ u \O}{f^{p+1}}) := u$ is a bijection. In the proof of \cite[Theorem 10.6]{BatCox} it is shown that the subspace
\be
H^0(\BP, \hat \O_{\BP}^m (p X)) + d H^0(\BP, \hat \O_{\BP}^{m- 1}(p X)) \subset H^0(\BP, \hat  \O_{\BP}^m ((p+1) X)) 
\ee
maps via $\psi$ to 
\be
J(f)_{p\b} \subset S_{p\b}
\ee
for any $p \geq 0$.

\begin{itemize}
	\item Suppose $p\geq m+1$. If \[ R(f)_{p\b }=   S_{p\beta }/ J(f)_{p\beta }= \frac{H^0(\BP,\hat  \Omega_\BP^m ((p+1)X))}{H^0(\BP,\hat  \Omega_\BP^m (pX)) + dH^0(\BP, \hat  \Omega_\BP^{m-1} (pX))} \] were non-zero, then $\frac{H^0(\BP,\hat  \Omega_\BP^m ((p+1)X))}{dH^0(\BP,\hat  \Omega_\BP^{m-1} (pX))}$ would be non-zero contributing non-trivially\footnote{The de Rham complex of $\BP$ with poles of arbitrary order along $X$ has the pole order filtration and this filtered complex is quasi-isomorphic to the de Rham complex of $\BP \setminus X$.
 } to $H^m(\BP\setminus X)$, which contradicts the isomorphism \eqref{bca}.
 
	\item Suppose $p=m$. Then \cite[Theorem 10.6]{BatCox} gives an isomorphism 
\be
  R(f)_{m\b}\cong \text{Gr}_F^0 H^m(\BP\setminus X).
\ee
By \cite[Corollary 10.12]{BatCox}, there is a natural isomorphism
\be
\text{Gr}_F^0 H^m(\BP\setminus X) \cong H_{\text{pr}}^{-1, m}(X) =0,
\ee
which implies that $R(f)_{m\b}=0$.
\end{itemize}

%{\color{blue}
%Since we assume that the map 
%\be
%H^{m-2}(\BP) \xrightarrow{\cup [X]} H^{m}(\BP)
%\ee
%is an isomorphism, we have an isomorphism
%\be
%H^m(\BP\setminus X) \simeq H_{\text{pr}}^{m-1}(X)
%\ee
%by the Gysin long exact sequence \cite[Proposition 10.11]{BatCox}. Thus, by comparing the isomorphism \eqref{bca} with the isomorphism \eqref{BCe} of \cite[Theorem 10.13]{BatCox}
%\be
%\bigoplus_{p=0}^m R(f)_{p\b} \simeq H^m(\BP\setminus X) \simeq H_{\text{pr}}^{m-1}(X) \simeq \bigoplus_{p=0}^{m-1}R(f)_{p\b},
%\ee
%we conclude that
%\be
%R(f)_{m\b}=0.
%\ee
%}
\end{proof}

\begin{corollary}\label{ebt}
 Under the same assumption with Theorem \ref{thm:BatCox}, there is a graded vector space isomorphism (depending on $\Omega$)
 \be
 \phi_\O: A(f) \xrightarrow{\cong} H^{m-1}_{\text{pr}}(X, \C).
 \ee
\end{corollary}
 
Now by Theorem \ref{thm:BatCox}, we have isomorphisms $\C\cong R(f)_0\cong H_{\text{pr}}^{m-1,0}(X)$ and $R(f)_{(m-1)\b}\cong \dim_\C H_{\text{pr}}^{0,m-1}(X)$. 
Since $\dim_\C H_{\text{pr}}^{m-1,0}(X) = \dim_\C H_{\text{pr}}^{0,m-1}(X)$ by the Hodge symmetry, we conclude that $\dim_\C R(f)_{(m-1)\b} =1$. 
Thus there is an isomorphism $R(f)_{(m-1)\b}\xrightarrow{\cong} \C$. We have the following analogue of Proposition \ref{prop:CY trace}.

%Using the isomorphism $\phi_\Omega$ in Corollary \ref{ebt} and the LG ring structure on $A(f)$, one can put another ring structure on $H^{m-1}_{\text{pr}}(X)$, which we call \textit{the LG (Landau--Ginzburg) ring structure}.
%by the celebrated theorem of Griffiths \cite{Gr69} (in the projective smooth hypersurface case) and Batyrev--Cox's theorem \cite[Theorem 10.13]{BatCox} for a quasi-smooth ample toric hypersurface in a complete simplicial toric variety. 
%More explicitly, one has an isomorphism between $A(f)$ and $H_{\text{pr}}^{m-1}(X)$; see \eqref{BCe} for details and the LG ring structure is induced from the ring structure of $A(f) \subset R(f)$. 
%Our main result is that these two ring structures are equivalent; see Theorem \ref{mthm} for a precise statement.
\begin{proposition}\label{mpro}
We further assume that $X$ is a non-degenerate hypersurface in $\BP$ (for example, see \cite[Definition 3.1]{Loyola}).
There exists an isomorphism $\Tr_{\text{LG}}:R(f)_{(m-1)\b}\xrightarrow{\cong} \C$ such that the following diagram
\[\begin{tikzcd}[row sep=0.5em, column sep=5em]
	{R(f)_{a\b} \times R(f)_{b\b}} & {R(f)_{(m-1)\b}} \\
	&& \mathbb{C} \\
	{H_{\text{pr}}^{m-1-a,a}(X) \times H_{\text{pr}}^{m-1-b,b}(X)} & {H_{\text{pr}}^{m-1,m-1}(X)}
	\arrow["{\text{mul}}", from=1-1, to=1-2]
	\arrow["{\phi_\O \times \phi_\O}", from=1-1, to=3-1, swap]
	\arrow["\Tr_{\text{LG}}", from=1-2, to=2-3]
	\arrow["\cup", from=3-1, to=3-2]
	\arrow["{\int_X}"', from=3-2, to=2-3]
\end{tikzcd}
\]
 commutes, where
\be
\text{mul}(u,v) = (-1)^b u \cdot v \qquad \text{for}\quad  u \in R(f)_{a\beta},\ v \in R(f)_{b\beta}.
\ee
% up to sign, i.e.,
%\be
%Tr(u \cdot v)= (-1)^{b} \int_X \phi(u) \wedge \phi(v), \quad u \in R(f)_{ar}, v \in R(f)_{br}
%\ee
whenever $a+b =m-1$.\footnote{Note that the bilinear map $\text{mul}$ is symmetric (respectively, skew-symmetric) if $m-1$ is even (respectively, $m-1$ is odd).}
\end{proposition}
%By the Batyrev--Cox theorem \cite[Theorem 10.13]{BatCox}, $\C=R(f)_0$ is isomorphic to $H_{\text{pr}}^{m-1,0}(X)$ and $R(f)_{(m-1)\b}\simeq \dim_\C H_{\text{pr}}^{0,m-1}$. 
%Since $\dim_\C H_{\text{pr}}^{m-1,0} = \dim_\C H_{\text{pr}}^{0,m-1}$ by the Hodge symmetry, we conclude that $\dim_\C R(f)_{(m-1)\b} =1$. 
For the proof, we will need Villaflor's toric Carlson--Griffiths theorem (see Theorem \ref{CD}), so we defer the proof to  the end of Subsection \ref{ss3.2}; in particular, see \eqref{LGtr} for an explicit definition of $\Tr_{\text{LG}}$.

Now, motivated by Proposition \ref{mpro}, we define a bilinear pairing (analogous to  the CY definition \eqref{cypair} and Proposition \ref{prop:CY trace})
\be
R(f)_{a\b} \times R(f)_{b\b} \xrightarrow{\text{mul}} R(f)_{(m-1)\b} \xrightarrow{\Tr_{\text{LG}}} \C
\ee
by
\be
\langle u, v \rangle_{\text{LG}} := \Tr_{\text{LG}}(u \cdot v)
\ee
whenever $a+b =m-1$, which induces a symmetric bilinear pairing
\be
\langle - , - \rangle_{\text{LG}}: A(f) \times A(f) \to \C,
\ee
by declaring that $\langle u, v\rangle = 0$ unless $\deg u+ \deg v= (m-1)\b$.

\begin{definition}
We call the triple $(A(f),  \bullet_{\text{LG}} , \langle - , - \rangle_{\text{LG}})$ the \emph{LG Frobenius algebra} of $f$. By transport of structure via $\phi_\Omega$, it defines the \emph{LG Frobenius algebra} $(H^{m-1}_{ \text{pr} }(X),  \bullet_{\text{LG}} , \langle - , - \rangle_{\text{LG}})$.
\end{definition}

\section{Main theorem and application}\label{sec3}
%{ via the Barannikov--Kontsevich method}

The main goal of this section is to prove Theorem \ref{mto} that the CY Frobenius algebra structure on  $H_{\text{pr}}^{m-1}(X)$ is equivalent to the LG Frobenius algebra structure on $H_{\text{pr}}^{m-1}(X)$.
The main ingredient of the proof is to verify that Batyrev--Cox's map is a ring isomorphism between $H_{\text{pr}}^{m-1}(X)$ with CY ring structure and $A(f)$ with LG ring structure by using Villaflor's generalization \cite{Loyola} of the Carlson--Griffiths theorem \cite{CG}; we state the main result as Theorem \ref{mthm} in Subsection \ref{ss3.1}. In Subsection \ref{ss3.2}, we prove Theorem \ref{mthm} and Proposition \ref{mpro}. In Subsection \ref{ss3.3}, we deduce a formal non-trivial Frobenius manifold structure on $(A(f), \bullet_{\text{LG}}, \langle -, - \rangle_{\text{LG}})$ and establish Corollary \ref{mtt}. In Subsection \ref{ss3.4}, we discuss the case where $\BP=\BP^{r-1}$ and provide  Example \eqref{keyexample} where the polynomial $f$, which is viewed as a function $\C^r \to \C$, may have non-isolated and non-compact critical locus. Thus, our construction yields a Frobenius manifold in a setting not covered by K. Saito's work on isolated hypersurface singularities or M. Saito's generalization \cite{MS}.
%{\color{blue}
%This may have a criticism from a referee regarding the construction is genuinely new; it might be ``just a transport of Frobenius structure".
%}

%We show that Theorem \ref{mthm} provides an identification of the Frobenius algebra structure on $H_{\text{pr}}^{m-1}(X)$ with CY ring structure and the cup product and the Frobenius algebra structure on $A(f)$ with LG ring structure and the torc residue map.
%As its application, when $X_f$ is smooth, we are able to provide a formal Frobenius manifold structure on the ring $A(f)$ using the Barannikov--Kontsevich construction of formal Frobenius manifolds on $H_{\text{pr}}^{m-1}(X)$ in \cite{BK}. 

%{\color{blue}
%Moreover, we use the Khovanskii-Pukhlikov presentation of the cohomology ring $A(\BP)$ of $\BP$ to provide a formal Frobenius manifold structure on $A(f) \oplus \frac{A(\BP)}{\C\cdot H^{m+\delta}(\BP)}$ where $\delta=0$ (respectively, $\delta=1$) for even $m$ (respectively, odd $m$). 
%}

%\subsection{CY Frobenius algebras and Frobenius manifolds of Barannikov--Kontsevich}\label{ss3.1}
\subsection{Comparison of two Frobenius algebras}\label{ss3.1}

%According to the Hodge decomposition of the primitive middle-dimensional cohomology, we have a graded vector space isomorphism
%\be
%H_{\text{pr}}^{m-1}(X, \C) = \bigoplus_{a=0}^{\infty} H^a_{\text{pr}}(X, \hat \O_X^{m-1-a})\xrightarrow{c_{\O_X}^{-1}} \bigoplus_{a=0}^{\infty}H^a_{\text{pr}}(X, \wedge^a \cT_X)
%\ee
%and the ring structure \eqref{CYR} of $ \oplus_{a=0}^{\infty}H^a_{\text{pr}}(X, \wedge^a \cT_X)$ transports to the CY ring structure on $H_{\text{pr}}^{m-1}(X, \C)$.
A natural question is to compare the two (CY versus LG) ring structures on $H_{\text{pr}}^{m-1}(X, \C)$. Our main result is that these two ring structures on $H_{\text{pr}}^{m-1}(X, \C)$ are in fact isomorphic.
%for even $n$ and anti-equivalent for odd $n$.
\begin{theorem}\label{mthm}
Let $X\subset \BP_\Sigma $ be a smooth\footnote{We assume smoothness because the result of Barannikov and Kontsevich is known only for a Calabi--Yau manifold.} Calabi--Yau hypersurface such that  \eqref{extraisom} is an isomorphism. For a given element $\O \in H^0(\BP_\Sigma, \hat \Omega_{ \BP_\Sigma }^m(\beta)  ) $, there is a nowhere-vanishing holomorphic volume form $\O_X$ on $X$ such that the vector space isomorphism
\be
\Phi:=c_{(-1)^{m-1}\O_X}^{-1} \circ \phi_\O: A(f) \xrightarrow{\cong}H_{\text{pr}}^0 (\PV(X))
\ee
is a ring isomorphism, i.e.,
\be
\Phi([g] \cdot [h]) = \Phi([g]) \cdot \Phi([h])
\ee
where $g,h\in A_\b$ and $[- ]$ is the equivalence class modulo $J(f) \cap A_\b$.
\end{theorem}

\begin{note}
We have the following diagram:
\[\begin{tikzcd}
	& R(f):= {S/J(f)} &&& {H^\bullet_{\text{pr}}(\PV(X))} \\
	{} & {A(f)=\bigoplus_{a=0}^{\infty} R(f)_{a\b}} & {H^{m-1}_{\text{pr}}(X) } && {\bigoplus_{a=0}^\infty H^a_{\text{pr}}(X, \wedge^a \cT_X)=H^0_{\text{pr}}(\PV(X)).}
	\arrow[hook', from=2-2, to=1-2]
	\arrow["\phi_\O", from=2-2, to=2-3]
	\arrow["{c_{(-1)^{m-1}\O_X}^{-1}}", from=2-3, to=2-5]
	\arrow[hook', from=2-5, to=1-5]
\end{tikzcd}\]	
\end{note}

We prove Theorem \ref{mthm} in Subsection \ref{ss3.2}.

\begin{corollary}\label{cfr}
The map $\Phi$ in Theorem \eqref{mthm} induces an isomorphism between the LG Frobenius algebra $(A(f),  \bullet_{\text{LG}} , \langle - , - \rangle_{\text{LG}} )$ and the CY Frobenius algebra $(H_{\text{pr}}^0 (\PV(X)), \bullet_{\text{CY}} , \langle - , - \rangle_{\text{CY}})$.
\end{corollary}

\begin{proof}
By Proposition \ref{prop:CY trace} and Proposition \ref{mpro} with $n=m-1$, we have the following commutative diagram
\[
\begin{tikzcd}[row sep=2em, column sep=5em]
    H^{a}_{\text{pr}}(X, \wedge^a \cT_X ) \times H^{b}_{\text{pr}}(X, \wedge^b \cT_X )\arrow{r}{\wedge^*} \arrow{d}[swap]{c_{\O_X} \times c_{\O_X}} & H^{m-1}_{\text{pr}}(X, \wedge^{m-1} \cT_X )\arrow{dr}{\Tr_{\text{CY}}} & \\
    {H_{\text{pr}}^{m-1-a,a}(X) \times H_{\text{pr}}^{m-1-b,b}(X)} \arrow{r}{\cup} & {H_{\text{pr}}^{m-1,m-1}(X)} \arrow{r}{\int_X} & \mathbb{C} \\
    {R(f)_{a\b} \times R(f)_{b\b}} \arrow{r}{\text{mul}} \arrow{u}{\phi_\O \times \phi_\O}  & {R(f)_{(m-1)\b}} \arrow{ur}[swap]{\Tr_{\text{LG}}} &
\end{tikzcd}
\]
for $a+b=m-1$. Replacing $\Omega_X$ by $(-1)^{m-1}\Omega_X$ doesn't change the pairing $\langle - , - \rangle_{\text{CY}}$. Combined with Theorem \ref{mthm},  the result follows.
\end{proof}

By definition of LG and CY Frobenius algebra structures on $H^{m-1}_{\text{pr}}(X)$, Theorem \ref{mto} immediately follows.
%On the other hand, the ring $A(f)$ naturally embeds into $R(f)=S/J(f)$ as rings (the LG ring structure). 

\begin{remark}\label{FTYsign}
During the preparation of this paper, we were informed that \cite[Theorem 3.5, Theorem 3.6]{FTY} presents a similar result to Theorem \ref{mthm} and Corollary \ref{cfr}, specifically in the case where $\BP=\BP^{r-1}$ is the projective space and $X_f$ is a smooth hypersurface in $\BP$. However, the assumption of smoothness for $X_f \subset \BP^{r-1}$ implies that $f$ has an isolated singularity at the origin, making it essential to consider non-projective cases of $\BP$ in order to discuss non-isolated singularities. Furthermore, our map $\Phi$ in Theorem \ref{mthm} is a ring isomorphism, resolving a sign ambiguity in the map $r'$ in \cite[Theorem 3.5]{FTY} even in the case of $\BP^{r-1}$.
\end{remark}

%\subsection{Toric Carlson--Griffiths theorem of Loyola}\label{ss2.3}
\subsection{Proofs of main results} \label{ss3.2}
%\be
%-K_{\BP}=\sum_{\rho \in \Sigma(1)} D_\rho
%\ee
Here we prove Theorem \ref{mthm} and Proposition \ref{mpro} by expanding upon the toric Carlson--Griffiths theorem due to Villaflor Loyola \cite{Loyola}.

Since $X_f$ is a quasi-smooth hypersurface in $\BP_\Sigma$, the collection $\cU=\{U_i : i=0,\cdots, r\}$  of open sets, where $U_i=\{ \ud a \in X_f  \mid   f_i(\ud a) \neq 0 \}$, is a contractible open covering of $X_f$, called the Jacobian covering of $X_f$. If $J=(j_0, j_1, \cdots, j_p)$ is an index set with size $|J|=p+1$, then we set
\be
U_J= U_{j_0} \cap \cdots \cap U_{j_p}.
\ee
Given a vector field $Z$ on $\C^{r}$, let $\iota(Z)$ denote the operation of contraction with $Z$.
Given a multi-index $J=(j_0,j_1, \cdots, j_p)$, let
\be
\O_J:= \iota\left (\frac{\partial}{\partial z_{j_p}}\right) \cdots  \iota\left (\frac{\partial}{\partial z_{j_0}} \right ) \O \qquad \text{and} \qquad f_J=f_{j_0}\cdots f_{j_p}.
\ee

Villaflor Loyola described the residue map in terms of \v{C}ech cocycle with respect to the Jacobian covering, which generalizes Carlson--Griffiths' result for smooth projective hypersurfaces in \cite[Proposition, page 7]{CG}:

\begin{theorem}\label{CD}\cite[Theorem 8.1]{Loyola} (Toric Carlson--Griffiths Theorem in the Calabi--Yau case)
Let $\BP_\Sigma$ be an $m$-dimensional projective simplicial toric variety with anti-canonical class $\b \in \Cl(\Sigma)$.
Let $X_f$ be a quasi-smooth ample hypersurface in $\BP_\Sigma $  of degree $\b=\deg (f) \in \Cl(\BP_\Sigma)$.
For $p \in \{0,1, \cdots, m-1\}$ and $F(\ud z) \in S_{p\b}$, one has\footnote{When $\BP$ is the projective space, note that Carlson--Griffiths' formula in \cite[Proposition in Section 3.b]{CG} has an additional sign factor $(-1)^{\frac{p(p+1)}{2}}$ compared to Villaflor's formula. This leads to a sign ambiguity of the map $r'$ in \cite[Theorem 3.5]{FTY} mentioned in Remark \ref{FTYsign}. 
}
\be
\res \left(\frac{F(\ud z) \Omega}{f^{p+1}}\right) =\frac{(-1)^{m-1}}{p! } \left\{ \frac{F(\ud z) \O_J}{f_J}\right\}_{|J|=p+1} \in H^p (\cU, \hat\Omega_{X_f}^{m-1-p}).
\ee
%where $\O_J:= \iota_{\frac{\partial}{\partial z_{j_p}}} (\cdots \iota_{\frac{\partial}{\partial z_{j_0}}} (\O)\cdots ), f_J:=f_{j_0} \cdots f_{j_p}$ and $\cU=\{U_i\}_{i=1}^r$ is the Jacobian covering restricted to $X_f$, given by $U_i =\{ f_i \neq 0\} \cap X_f\subset \BP_\Sigma$.
\end{theorem}

%\begin{proposition}\cite[Proposition, page 7]{CG}\label{CD} We have
%\be
%\mathrm{res} \bigg(\bigg[\frac{(-1)^a a ! g(\ud x) \Omega}{f(\ud x)^{a+1}}\bigg]\bigg)
%=\Big\{ (-1)^{n+a+ \frac{a(a+1)}{2}}\frac{g \Omega_J}{f_J} \Big\}_{|J|=a} \quad \text{in} \quad H^a(X, \wedge^{n-a} \cT_X^*).
%\ee
%\end{proposition}

Now we specify a Calabi--Yau volume form $\O_X$ on $X=X_f(\C)$ using $\O$:
\begin{lemma}\label{CYform}
For each $i, j$, we have
\be
\frac{\frac{\partial}{\partial z_i} \vdash \Omega}{ f_i} = \frac{\frac{\partial}{\partial z_j} \vdash \Omega}{ f_j }  \quad \text{on} \quad U_{(i,j)}=U_i \cap U_j.
\ee
This implies that $\Omega_{X}|_{U_j} :=\frac{\frac{\partial}{\partial z_j} \vdash \Omega}{ f_j}= \frac{\O_{j}}{f_{j}} $ for each $j$ glues together to define a nowhere-vanishing holomorphic $n$-form $\Omega_X$ on the Calabi--Yau manifold $X$. 
In other words, the volume form $\O_X$ is given by $\res\left ((-1)^{m-1} \frac{ \O}{f}\right )$.

%This implies that $\Omega_{X_f}:=\frac{\frac{\partial}{\partial x_i} \vdash \Omega}{ f_i}$ on $U_i$ for each $i$ glues together to define a non-vanishing holomorphic $n$-form on the Calabi--Yau manifold $X_f$. 
\end{lemma}
\begin{proof}
By \cite[Corollary 7.2]{Loyola}, we have (see \cite[Corollary 7.1]{Loyola} for the definition of $ V^\b_j$)
\be
\O_{j} \wedge df + (-1)^m f_{j} \O = (-1)^{(m-1)r} \cdot f \cdot V_{j}^\b.
\ee
Therefore we have
\be
\frac{\O_{j}}{f_{j}}\wedge \frac{df}{f} -(-1)^{(m-1)r} \cdot \frac{V_{j}^\b}{f_j}= \frac{(-1)^{m-1} \O}{f}
\ee
for any $j$ and hence we get
\be
\res\left(  \frac{\O_{j}}{f_{j}}\wedge \frac{df}{f} -(-1)^{(m-1)r} \cdot \frac{V_{j}^\b}{f_j} \right) = \res\left(  \frac{(-1)^{m-1} \O}{f}\right) =  \frac{\O_{j}}{f_{j}}
\ee
where the latter equality follows from Theorem \ref{CD}. On the other hand, because
\be
\frac{\O_{i}}{f_{i}}\wedge \frac{df}{f} -(-1)^{(m-1)r} \cdot \frac{V_{i}^\b}{f_i} =  \frac{(-1)^{m-1} \O}{f} = \frac{\O_{j}}{f_{j}}\wedge \frac{df}{f} -(-1)^{(m-1)r} \cdot \frac{V_{j}^\b}{f_j}
\ee
for each $i,j$, the claim follows. 
\end{proof}

\begin{lemma}\label{WD}
For $J$ with $|J|=p+1$ and each $0\leq k, l\leq p$, we have
\be \label{polyv}
 \frac{(-1)^k f_{j_k} \frac{\partial}{\partial z_{j_p}} \cdots \widehat{ \frac{\partial}{\partial z_{j_k}}} \cdots \frac{\partial}{\partial z_{j_0}} }{f_J}
 = \frac{(-1)^l f_{j_l} \frac{\partial}{\partial z_{j_p}} \cdots \widehat{ \frac{\partial}{\partial z_{j_l}}} \cdots \frac{\partial}{\partial z_{j_0}} }{f_J}
 \quad \text{on} \quad U_J.
\ee
\end{lemma}

\begin{proof}
Lemma \ref{CYform} implies that
\be
\frac{(-1)^k f_{j_k} \left( \frac{\partial}{\partial z_{j_p}} \cdots \widehat{ \frac{\partial}{\partial z_{j_k}}} \cdots \frac{\partial}{\partial z_{j_0}} \right)\vdash \O_{X}}{f_J}
=\frac{(-1)^l f_{j_l} \left( \frac{\partial}{\partial z_{j_p}} \cdots \widehat{ \frac{\partial}{\partial z_{j_l}}} \cdots \frac{\partial}{\partial z_{j_0}} \right)\vdash \O_{X}}{f_J}
%= \frac{(-1)^jf_{j_j} \frac{\partial}{\partial x_{j_p}} \cdots \hat{ \frac{\partial}{\partial x_{j_j}}} \cdots \frac{\partial}{\partial x_{j_0}} }{f_J}
\quad \text{on} \quad U_J,
\ee 
which gives us the desired equality.
\end{proof}

\begin{proof}[Proof of Theorem \ref{mthm}]\label{proof of mtm}
In order to show that $\Phi = c_{(-1)^{m-1}\O_X}^{-1} \circ \phi_\O$ is a ring isomorphism, one first notes that the map $\Phi$ is given, using the \v{C}ech description in Theorem \ref{CD} of $\bigoplus_{a=0}^\infty H^a_{\text{pr}}(X, \wedge^a \cT_X)$ in terms of the Jacobian covering, as follows:
\bea \label{describe}
\Phi([g(\ud z)]) = \left\{ (-1)^{a}g(\ud z) \frac{P_{J_i}}{f_{J}}  \right\}_{|J|=a+1} \qquad \text{for}\quad  g(\ud z) \in S_{a\b}
\eea
where $P_{J_i}=(-1)^if_{j_i} \frac{\partial}{\partial z_{j_a}} \cdots \widehat{ \frac{\partial}{\partial z_{j_i}}} \cdots \frac{\partial}{\partial z_{j_0}} \in \G(U_J,\wedge^a \cT_X)$ and $f_{J}=f_{j_0}\cdots  f_{j_a}$ with the index set $J=(j_0, j_1, \cdots, j_a)$.
By Lemma \ref{WD}, we have
\be
 \frac{P_{J_k}}{f_{J}}   =  \frac{P_{J_l}}{f_{J}}  \quad \text{on} \quad U_J \subset X
\ee
for any $k,l \in \{0, \cdots, a\}$, which implies that $\Phi$ in \eqref{describe} is well-defined on $U_J$.
When $J=(j_0, j_1, \cdots, j_a)$, let $\tilde J$ denote $(j_1, \cdots, j_a)$.
Note that $\Phi([g(\ud z)])$ can be also written as
%\be
%\Phi([g(\ud x)]) = \{ (-1)^{a}g(\ud x) \frac{\tilde P_{J_0}}{f_{\tilde J}}  \}_{|J|=a+1}, \quad  g(\ud x) \in S_{a\b},
%\ee
%where $\tilde J=(j_1, \cdots, j_a)$ and $\tilde P_{J_0}=\frac{\partial}{\partial x_{j_a}} \cdots \frac{\partial}{\partial x_{j_1}}$.
\be
\Phi([g(\ud z)]) = \left\{ (-1)^{a}g(\ud z) \frac{\frac{\partial}{\partial z_{j_a}} \cdots \frac{\partial}{\partial z_{j_1}}    }{f_{j_1} \cdots f_{j_a}}  \right\}_{|J|=a+1}\qquad \text{for}\quad g(\ud z) \in S_{a\b}.
\ee

As in \cite[page 13]{CG}, let ``$\wedge$'' denote the natural product
\be
\cC^a(\cU|_X, \wedge^b\cT_X) \times \cC^c(\cU|_X, \wedge^d \cT_X)  \to \cC^{a+c}(\cU|_X, \wedge^{b+d}\cT_X) 
\ee
given by the ``front $a$-face, back $c$-face'', followed by exterior multiplication of polyvector fields. Then the twisted product is given by
\bea \label{sign}
\a^b_a \cdot  \a^d_c = (-1)^{ad}\a^b_a \wedge \a^d_c,
\eea
which represents the cup product on the level of hypercohomology. Since
\be
\Phi([g(\ud z)]) &=& (-1)^{a}\left \{ g(\ud z) \frac{P_{J_a}}{f_J} \right \}_{J=(j_0, \cdots, j_a)},  \quad  g(\ud z) \in S_{a\b},\\
\Phi([h(\ud z)]) &=& (-1)^{b}\left \{ h(\ud z) \frac{P_{J'_0}}{f_{J'}} \right \}_{J'=(j'_0, \cdots, j'_b)},  \quad  h(\ud z) \in S_{b\b},
\ee
the cup product of $\Phi([g(\ud z)]) \in \cC^a(\cU|_X, \wedge^a\cT_X) $ and $\Phi([h(\ud z)])\in \cC^b(\cU|_X, \wedge^b\cT_X) $ is given by
\be
&& \Phi([g(\ud z)]) \cdot \Phi([h(\ud z)])  \\
&\stackrel{\eqref{sign}}{=}&(-1)^{ab}\left \{(-1)^{a+b}
g(\ud z) \frac{(-1)^a f_{s} \frac{\partial}{\partial z_{j_{a-1}}} \cdots \frac{\partial}{\partial z_{j_0}}}   {f_{s} f_{j_{a-1}}\cdots f_{j_0}}  
h(\ud z) \frac{(-1)^0 f_{s} \frac{\partial}{\partial z_{j_b'}} \cdots \frac{\partial}{\partial z_{j_1'}} }   {f_{j_b'} \cdots f_{j_1'} f_{s}}  \right \}_{(j_0,\cdots , j_{a-1} , s , j'_1, \cdots, j'_b)} \\
&=&\left \{(-1)^{a+b}
g(\ud z) h(\ud z)
\frac{(-1)^a f_{i_a} \frac{\partial}{\partial z_{i_{a+b}}} \cdots \frac{\partial}{\partial z_{i_{a+1}}} \cdot \frac{\partial}{\partial z_{i_{a-1}}} \cdots \frac{\partial}{\partial z_{i_0}}  }   
{f_{i_{a+b}} f_{i_{a+b-1}}\cdots f_{i_1} f_{i_0}}  
\right \}_{I=(i_0, \cdots,i_a,i_{a+1},\cdots, i_{a+b})}.
\ee
%Then the topological cup-product of $\Phi([g(\ud x)])$ and $\Phi([h(\ud x)])$, where $\deg g= a(n+2)$ and $\deg h=b(n+2)$, is given by
%\be
%&&
%\Phi([g(\ud x)]) \cdot \Phi([h(\ud x)])=\Big\{(-1)^{ab}(-1)^{n+a+\frac{a(a+1)}{2}}g(\ud x) \frac{\tilde P_{J_0}}{f_{\tilde J}}   (-1)^{n+b+\frac{b(b+1)}{2}}h(\ud x) \frac{\tilde P_{J'_0}}{f_{\tilde J'}}\Big\}_{\tilde J' \sqcup \tilde J} \\
%%&=&
%\ee
%where $\tilde J'\sqcup \tilde J = (j_1', \cdots, j_b', j_1, \cdots, j_a), \tilde P_{J_0}=\frac{\partial}{\partial x_{j_a}} \cdots \frac{\partial}{\partial x_{j_1}}$ and $\tilde P_{J_0'}=\frac{\partial}{\partial x_{j'_b}} \cdots \frac{\partial}{\partial x_{j'_1}}$.
On the other hand, we have
\be
\Phi([g(\ud z)]\cdot[ h(\ud z)]) = \left \{ (-1)^{a+b}g(\ud z)h(\ud z) \frac{P_{J_a}}{f_{J}}\right\}_{|J|=a+b+1}.
\ee
Therefore we have $\Phi([g] \cdot [h]) =  \Phi([g]) \cdot \Phi([h])$, which finishes the proof of Theorem \ref{mthm}.
\end{proof}

\vspace{1em}

\begin{proof}[Proof of Proposition \ref{mpro}]
The proof amounts to careful application of results of Villaflor \cite{Loyola}.

For $f \in S = \C [ \ud z  ] $, we define an ideal
 $J_0(f):=\langle z_1 f_1, \cdots, z_m f_m, \cdots, z_r f_r \rangle$ of $S$ and set $R_0(f):=S/J_0(f)$.

Let $\rho_{i_1}, \cdots, \rho_{i_m}$ be linearly independent primitive generators of the rays of $\Sigma$. Then for $I=\{i_1, \cdots, i_m\}$, \cite[Section 6]{Loyola}
%\cite[Definition 6.1, Definition 6.2, Theorem 6.1, and Definition 6.3]{Loyola}
gives a definition of the toric Hessian of $f$ as
\bea \label{wd}
\operatorname{Hess}_\Sigma(f):= \frac{\operatorname{Hess}_\Sigma^I (f)}{\det (\rho_I)} \in R_0(f)_{m\beta},
\eea
which turns out to be independent of the choice of $I$.
Since $\dim_\C R_0(f)_{m\b}=1$  and $\operatorname{Hess}_\Sigma(f)$ is non-zero by \cite[Proposition 3.3, Corollary 6.1]{Loyola}, for $V  \in R_0(f)_{m \b}$, there is a unique number $c_V \in \C$ such that
\be
V \equiv c_V \cdot \operatorname{Hess}_\Sigma(V ) \mod J_0(f).
\ee

Motivated by \cite[Corollary 6.1]{Loyola}, we define the LG trace pairing $\Tr_{\text{LG}}: R(f)_{(m-1)\b} \to \C$ by
\bea \label{LGtr}
\Tr_{\text{LG}} (U ) :=-(2 \pi i)^{m-1} (-1)^{\frac{m(m-1)}{2}} c_{z_1 \cdots z_r U } \cdot m ! \cdot \operatorname{Vol}(\Delta),
\eea
where $\Delta$ is the convex polyhedron associated to the anti-canonical divisor $\b$. Here we used that if $U  \in  S_{(m-1)\b}$, then $z_1 \cdots z_r U \in S_{m\b}$ and that $z_1\cdots z_r J(f)_{(m-1)\beta } \subset J_0(f)_{m\beta }$ holds.  

Now let us calculate the cup product in \v{C}ech cohomology
\be
\cC^a(\cU|_X, \hat\O^{b}_X) \times \cC^b(\cU|_X, \hat\O^{a}_X)  \to \cC^{m-1}(\cU|_X, \hat\O^{m-1}_X), \quad a+b =m-1.
\ee
By Theorem \ref{CD}, we have
 \be
 \phi_\O ([u(\ud z)]) 
 := \res \left((-1)^a a!  \frac{u(\ud z) \Omega}{ f^{a+1}} \right) 
 = (-1)^{m-1+a} \left\{ \frac{u(\ud z) \O_J}{f_J}\right\}_{|J|=a+1} \in H^a (\cU, \hat\Omega_{X}^{m-1-a})=H^a (\cU, \hat\Omega_{X}^{b})
\ee
for $ u(\ud z) \in R(f)_{a\b}$. Then for $ u(\ud z) \in R(f)_{a\b}$ and  $v(\ud z) \in R(f)_{b\b}$ (with $a+b=m-1$), the cup product is given by 
\be
  \phi_\O ([u(\ud z)])  \cup  \phi_\O ([v(\ud z)]) 
  &=& (-1)^{a^2+a+b} \left\{ \frac{u(\ud z) v(\ud z) \Omega_{Rs}\wedge \Omega_{sT}}{f_R f_s^2 f_T}\right\}_{|L|=m}  \\
   &=&(-1)^{b} \left\{ \frac{u(\ud z) v(\ud z) \Omega_{Rs}\wedge \Omega_{sT}}{f_R f_s^2 f_T}\right\}_{|L|=m}
   \in \cC^{m-1}(\cU|_X, \O_X^{m-1}) 
\ee
where the sign $(-1)^{a^2}$ arises from \eqref{sign} and the multi-index $L$ is partitioned to
\be
L=(r_0, \cdots, r_{a-1}, s, t_1, \cdots, t_b), \quad R=(r_0, \cdots, r_{a-1}), \quad  T=(t_1, \cdots, t_b).
\ee
 Note that $\dim_\C R(f)_{(m-1)\b}=1$ and $u(\ud z)v(\ud z) \in R(f)_{(m-1)\b}$.
Since one can show $\O_{Rs} \wedge \O_{sT} = \det (\rho_L) \widehat{z_L} \O_s$ from \eqref{vf}, we have
\be
\left\{ \frac{u(\ud z) v(\ud z) \Omega_{Rs}\wedge \Omega_{sT}}{f_R f_s^2 f_T}\right\}_{|L|=m}
=\left\{ \frac{u(\ud z) v(\ud z)  \det (\rho_L) \widehat{z_L} \O_s}{f_L f_s}\right\}_{|L|=m}.
\ee
We recall the reside map $\res$ of the Poincar\'e regular sequence in \cite[Definition 10.1]{Loyola} and the coboundary map \cite[(19)]{Loyola}
\be
\tau: H^{m-1}(X, \hat \O_{X}^{m-1}) \xrightarrow{\cong} H^m(\BP, \hat\O_{\BP}^m).
\ee
For simplicity of notation, we write $L=(\ell_1, \cdots, \ell_m)$ and $\tilde L= (\ell_0, \ell_1, \cdots, \ell_m)$ below.
In order to compute $\tau$, we consider a lift of the cocycle in the Poincar\'e regular sequence
\be
\eta_L:=\frac{u(\ud z) v(\ud z)  \det (\rho_L) \widehat{z_L} \O_s}{f_L f_s} \wedge \frac{df}{f}
-\frac{u(\ud z) v(\ud z)  \det (\rho_L) \widehat{z_L} V_s^\b}{f_L f_s}, 
\ee
where we refer to \cite[Corollary 7.1]{Loyola} for the notation $V^\b$ and \cite[Corollary 7.2]{Loyola} for the underlying idea;  it gives $\res (\eta_L) = \left\{ \frac{u(\ud z) v(\ud z)  \det (\rho_L) \widehat{z_L} \O_s}{f_L f_s}\right\}_{|L|=m}$. Using this we can compute the coboundary map (as in the proof of \cite[Proposition 11.1]{Loyola}) by \cite[Corollary 7.2, (10)]{Loyola} 
{\begin{small}
  \begin{align}\label{oned}
 \star & = \tau \left( \left\{ \frac{u(\ud z) v(\ud z)  \det (\rho_L) \widehat{z_L} \O_s}{f_L f_s}\right\}_{|L|=m} \right)= \left\{  \frac{u(\ud z)v(\ud z) \widehat{z_{\tilde L}} \sum_{i=0}^m (-1)^i \det (\rho_{\tilde L \setminus \{\ell_i\}}) z_{\ell_i} f_{\ell_i}}{f_{\tilde L} f} \O\right\}_{|\tilde L|=m+1} \in H^{m}(\cU, \hat \O_{\BP}^m).
\end{align}
\end{small}}
%{\begin{small}
%\bea \label{oned}
%\qquad \qquad  \star= \tau \left( \left\{ \frac{u(\ud z) v(\ud z)  \det (\rho_L) \widehat{z_L} \O_s}{f_L f_s}\right\}_{|L|=m} \right)
%= \left\{  \frac{u(\ud z)v(\ud z) \widehat{z_{\tilde L}} \sum_{i=0}^m (-1)^i \det (\rho_{\tilde L \setminus \{\ell_i\}}) z_{\ell_i} f_{\ell_i}}{f_{\tilde L} f} \right\}_{|\tilde L|=m+1} 
%\eea
%\end{small}}
Now without loss of generality one can take $I = \{1, \cdots, m\}$. Then let $\cU^I$ be the covering associated to $V(f, z_1f_1, \cdots, z_m f_m)= \phi \subseteq \BP$, i.e., given by $\cU_0=\{ f\neq 0\}$ and $\cU_i =\{ z_i f_i \neq 0\}$ for $i=1, \cdots, m$. This covering makes sense because of the result of Cox \cite[Proposition 5.3]{Cox} (see \cite[Proposition 3.2]{Loyola} for the form we use), since $X$ is assumed to be a non-degenerate ample hypersurface in $\BP$. Moreover, \cite[Proposition 3.2]{Loyola} says that
%the ideal $J_0(f):=\langle z_1 f_1, \cdots, z_m f_m, \cdots, z_r f_r \rangle$ is given by
\be
J_0(f) = \langle f, z_1 f_1, \cdots, z_m f_m \rangle.
\ee
Now let us write \eqref{oned} relative to the covering $\cU^I$ using the non-trivial Euler relations; we follow \cite[page 115]{Loyola} to compute as follows
\be
 \star=\left\{  \frac{u(\ud z)v(\ud z) \widehat{z_{1, \ldots, m}} \det (\rho_1 | \cdots | \rho_m)}{f \cdot  f_1 \cdots f_m}\Omega \right\}_{\{0\} \cup I} 
 = \left\{  \frac{u(\ud z)v(\ud z) z_1 \cdots z_r \det (\rho_1 | \cdots | \rho_m)}{f \cdot  z_1f_1 \cdots z_mf_m}\Omega \right\}_{\{0\} \cup I} 
 \in H^m (\cU^I, \hat \O_{\BP}^m).
\ee
By the definition \eqref{LGtr}, we have
\be
\Tr_{\text{LG}} ((-1)^b u \cdot  v ) = -(2 \pi i)^{m-1} (-1)^{\frac{m(m-1)}{2}}  c_{(-1)^b z_1\cdots z_r u   v } m! \operatorname{Vol}(\Delta).
\ee
On the other hand, we have the equality 
\be
-2 \pi i  \int _X \omega  =  \int_{\BP} \tau(\omega) =  (2\pi i)^m \Tr_{\BP}( \tau(\omega)) 
\ee 
for $\omega \in H^{m-1,m-1}(X)$ by \cite[Proposition 10.1]{Loyola}, where we use $\Tr_{\BP}$ given in \cite[Section 6]{Loyola}. From this we obtain
\begin{align*}
\int_X  \phi_\Omega(u)\cup \phi_\Omega(v) & =   - (2 \pi i)^{m-1} \Tr_{\BP } ( \tau( \phi_\Omega(u)\cup \phi_\Omega(v)) )\\
& =- (2 \pi i)^{m-1}\Tr_{\BP} \left((-1)^b \left\{  \frac{u(\ud z)v(\ud z) z_1 \cdots z_r \det (\rho_1 | \cdots | \rho_m)}{f \cdot  z_1f_1 \cdots z_mf_m}\O \right\}_{\{0\} \cup I} \right)\\
& = -(2 \pi i)^{m-1} (-1)^{\frac{m(m-1)}{2}} c_{(-1)^b z_1\cdots z_r u  v } m! \operatorname{Vol}(\Delta).
\end{align*}
where the last  equality follows from \cite[Corollary 6.1]{Loyola}. 
%\be
%\Tr \left((-1)^b \left\{  \frac{u(\ud z)v(\ud z) z_1 \cdots z_r \det (\rho_1 | \cdots | \rho_m)}{f \cdot  z_1f_1 \cdots z_mf_m}\right\}_{\{0\} \cup I}  \right)
%=-\frac{1}{(2 \pi i)^m} c_{(-1)^b z_1\cdots z_r u(\ud z) v(\ud z)} m! \operatorname{Vol}(\Delta).
%\ee
Note that we use the fact that \eqref{wd} is independent of $I$ and $\det (\rho_1 | \cdots | \rho_m) = \det(\rho_I)$.
Thus we conclude that $\Tr_{\text{LG}} ((-1)^b u \cdot  v )  = \int_X  \phi_\Omega(u)\cup \phi_\Omega(v)$.
%In the proof of \cite[Proposition 11.1]{Loyola}, the computation of $\int_X$ of the de Rham cohomology class of $ \left\{ \frac{F(\ud z) \Omega_{Rs}\wedge \Omega_{sT}}{f_R f_s^2 f_T}\right\}_{|L|=m}$
%for $F(\ud z) \in R(f)_{(m-1)\b}$ were given, from which we deduce the existence of $\Tr_{\text{LG}}:R(f)_{(m-1)\b}\xrightarrow{\cong} \C$ such that the diagram in Proposition \ref{mpro} commutes.
%Moreover,
%\be
%\O_{Rs}\wedge \O_{sT}= \det(\rho_L) \widehat{z_L} \O_{s}.
%\ee
\end{proof}

\subsection {Transport of formal Frobenius manifold structures}\label{ss3.3}

We briefly recall the definition of Frobenius manifolds and formal Frobenius manifolds. 
%We will only consider a pure even manifold.
\begin{definition}[Frobenius manifolds]
A \emph{Frobenius manifold} is a tuple $(M, \circ, e, E, g)$ where $M$ is a complex connected manifold with metric $g$, $\circ$ is a commutative and associative $\CMcal{O}_M$-bilinear multiplication $\cT_M \times \cT_M \to \cT_M$, and $e$ is a global unit vector field with respect to $\circ$, 
%and $E$ is another global vector field (called the Euler vector field) 
subject to the following conditions:
\begin{enumerate}[(1)]
\item (invariance) $g(X\circ Y, Z) = g(X, Y\circ Z)$,
\item (potentiality) the $(3,1)$-tensor $\nabla^g \circ$ is symmetric where $\nabla^g$ is the Levi-Civita connection of $g$,
\item the metric $g$ is flat, i.e. $\nabla^g$ is a flat connection, $[\nabla^g_X,\nabla^g_Y]=\nabla^g_{[X,Y]}$,
\item (flat identity) $\nabla^g e=0$.
%\item $\Lie_E(\circ) = \circ$ and $\Lie_E(g) = D \cdot g$ for some $D \in \C$.
\end{enumerate}
\end{definition}
One can similarly define a formal version of Frobenius manifolds by considering the formal structure sheaf and the formal tangent bundle instead of the holomorphic structure sheaf and the holomorphic tangent bundle.

\begin{definition}[formal Frobenius manifolds]\label{Frob}
	Let $M$ be a complex connected manifold of finite dimension, and $\ud t := \{t^\a\}$ be formal coordinates on an open subset $V \subset M$. Choose sufficiently small open subset $U$ of $V$ so that we could assume $\cO(U)=\mathbb{C}[\![\ud t]\!]$. On the local coordinates $(U, \ud t)$, let us write
	\[
		\partial_\alpha \circ \partial_\beta := \sum_\gamma A_{\alpha\beta}^\gamma \partial_\gamma.
	\]
	where $\{\partial_\a:=\partial/\partial_{t^\a}\}$ is a basis of $\cT_M$ and $A_{\a\b}^\g \in \bC[\![\ud t]\!]$ is a formal power series representing the 3-tensor field. 	Let $g$ be a non-degenerate symmetric bilinear pairing on $\mathbb{C}$-vector space spanned by $\{\partial_\alpha\}$. Then one can extend $g$ to a symmetric $\mathbb{C}\llbracket  \underline{t}\rrbracket $-linear pairing. Let $g_{\alpha\beta}:=g(\partial_\alpha,\partial_\beta) \in \C \
	$. Then $(M,\circ,g)$ is called a \emph{formal Frobenius manifold} if the following conditions are satisfied:
	\begin{enumerate}[(D1)]
		\item\label{D1} (associativity)
		\[
			\sum_\rho A_{\alpha\beta}^\rho A_{\rho\gamma}^\delta = \sum_\rho A_{\beta\gamma}^\rho A_{\rho\alpha}^\delta.
		\]
		\item\label{D2} (commutativity)
		\[
			A_{\beta\alpha}^\gamma = A_{\alpha\beta}^\gamma.
		\]
		\item\label{D3} (invariance) If we set $A_{\alpha\beta\gamma} = \sum_\rho A_{\alpha\beta}^\rho g_{\rho\gamma}$, then
		\[
			A_{\alpha\beta\gamma} = A_{\beta\gamma\alpha}.
		\]
		\item\label{D4} (flat identity) A distinguished element $\partial_0$ is the identity with respect to $\circ$:
		\[
			A_{0\alpha}^\beta =\boldsymbol\delta_\alpha^\beta \quad (\text{where $\boldsymbol\delta_\alpha^\beta$ is the Kronecker delta}).
		\]
		\item\label{D5}(potentiality)
		\[
			\partial_\alpha A_{\beta\gamma}^\delta= \partial_\beta A_{\alpha\gamma}^\delta.
		\]
	\end{enumerate}
\end{definition}

Barannikov--Kontsevich constructed a formal non-trivial Frobenius manifold structure
% (\cite[Section 1]{BK} for its definition) 
on $H_{\text{pr}}^0(\PV(X))$ extending the CY Frobenius algebra $(H_{\text{pr}}^0 (\PV(X)), \wedge, \langle -, - \rangle_{\text{CY}})$.

\begin{theorem}\cite{BK} \label{BKt}
There exists a 3-tensor $A_{\a\b}^\g(\ud s) \in \C\llbracket \ud s\rrbracket$ and the metric $g_{\a\b} \in \C$ (which is given by $\langle -,- \rangle_{\text{CY}}$) in a formal flat coordinate system $\ud s =(s_1, \cdots, s_\mu)$ on $H_{\text{pr}}^0 (\PV(X))$, where $\mu= \dim H_{\text{pr}}^0 (\PV(X))$, which provides a formal Frobenius manifold structure on $H_{\text{pr}}^0(\PV(X))$ extending\footnote{The basis $\{\o_\a: \a \in I\}$ of $H_{\text{pr}}^0(\PV(X))$ corresponding to the flat coordinates $\ud s$ satisfies
\be
\o_\a \cdot \o_\b =\sum_{\g \in I }A_{\a\b}^\g (\ud 0) \o_\g, \quad \langle \o_\a, \o_\b \rangle_{\text{CY}} = g_{\a\b}, \quad \o_\a, \o_\b \in H^0 _{\text{pr}} (\PV(X)).
\ee
} the CY Frobenius algebra $(H_{\text{pr}}^0 (\PV(X)), \wedge, \langle -, - \rangle_{\text{CY}})$. 
%This construction of formal Frobenius manifolds restricts to $H_{\text{pr}}^0(\PV(X))$.
\end{theorem}

Now we can transport the formal non-trivial Frobenius manifold structure to the Frobenius algebra $(A(f), \bullet_{\text{LG}}, \langle -, - \rangle_{\text{LG}})$: Corollary \ref{cfr} combined with Theorem \ref{BKt} implies Corollary \ref{mtt}. 

\subsection{Comparison with Isolated Singularities}\label{ss3.4}

Let us discuss Frobenius manifold structures on the case when $\BP=\BP^{r-1}$ is the projective space and $X_f \subset \BP^{r-1}$ is smooth:

\begin{remark}
Note that there is no direct relationship between the Jacobian algebra $R(f)$ and the total cohomology group $H^\bullet(\PV(X))$, but there is a concrete ring isomorphism between their subalgebras, namely, between $A(f)$ and $H^0(\PV(X))$ as we saw in Theorem \ref{mthm}.
%$\bigoplus_{a=0}^{\infty}H^a(X,\wedge^a\cT^*_X)$, 
The CY ring structure was put to use in \cite{BK} to construct a formal Frobenius manifold structure on $H^\bullet(\PV(X))$ in the context of Calabi--Yau B-model of mirror symmetry. On the other hand, the polynomial ring structure (the LG ring structure) on $S=\C[\ud z]$ was used to construct a Frobenius manifold structure on $R(f)$ by the theory of primitive forms and the higher residue parings associated to a universal unfolding of an isolated singularity $f$ \cite{Saito, Saitoh}; this Frobenius manifold structure is the main player in the context of Landau--Ginzburg B-model of mirror symmetry.

We observe that the Frobenius manifold structure on $H^\bullet(\PV(X))$ based on the CY ring structure (see \cite{BK} and its generalization \cite{Bara}) restricts to $H^0(\PV(X))$ and the Frobenius manifold structure on $R(f)$ based on the LG ring structure (see \cite{Saito, Saitoh} and also \cite{LLS}, \cite[Section 3]{CSKY} for its concrete algorithm) also restricts to $A(f)$.  Therefore, there are two types (CY versus LG) of constructions of Frobenius manifolds on $H^{r-2}_{\text{pr}}(X)$, each of which is an important invariant in the context of the $B$-model CY and LG mirror symmetry, respectively.

The CY Frobenius manifold structure depends on a Maurer--Cartan solution of the relevant dGBV algebra and a good opposite filtration to the Hodge filtration on $H^\bullet(X)$ and the LG Frobenius manifold structure depends on a universal unfolding of $f$ and a good opposite filtration to the ``Hodge filtration'' on $R(f)$.
It seems to be a subtle question to compare these two constructions on $H_{\text{pr}}^{r-2}(X) \cong A(f)$ explicitly.
The LG construction is more explicit in that it deals directly with a polynomial $f$ (while the CY construction deals with holomorphic polyvector fields on a manifold $X$), and can be made into a concrete algorithm (see \cite{LLS}) based on the Gr\"{o}bner basis by using the Jacobian ideal membership problem for $f$.

Theorem \ref{mto} and Corollary \ref{mtt} provide a non-trivial Frobenius manifold structure on the LG Frobenius algebra simply by transporting the CY Frobenius manifold structure. Hence a natural open question would be whether there exists an LG type construction meaning ``the theory of primitive forms and higher residue pairing'' (in the sense of \cite{LLS}) or the $L^2$-Hodge theoretic construction (in the sense of \cite{LW}) on $(A(f), \bullet_{\text{LG}}, \langle -, - \rangle_{\text{LG}} )$.
\end{remark}

Finally, we conclude with a concrete example of a function $f$ with a critical locus that is non-isolated and non-compact, yet to which our main theorem applies.

\begin{example}\label{keyexample}
We consider the projective $\BP^1$-bundle over $\BP^6$ associated to the vector bundle $\cO_{\BP^6}(2) \oplus \cO_{\BP^6}(3)$:
\be
\BP=\BP(\cO_{\BP^6}(2) \oplus \cO_{\BP^6}(3)),
\ee
which is a smooth projective Fano toric variety of dimension 7 over $\C$. In general, it is known that if $X$ is a smooth projective Fano variety and $D_1, \cdots, D_k$ are nef divisors on $X$ such that $-K_{X}-D_1 - \cdots - D_k$ is ample, then the projective bundle $\BP(\oplus_{i=1}^k \cO_X(D_i))$ is a Fano variety.
The homogeneous coordinate ring of $\BP$ is $\C[x_0, \cdots, x_6, y_1, y_2]=\C[z_1, \cdots, z_9]$, which is graded by the class group $\Cl(\BP)$ of $\BP$.
In this case, $\Cl(\BP)$ is isomorphic to $\Z^2$. In fact, we have
\be
\deg(x_i)=(1, 0) \quad \text{for}\quad i =0, \cdots, 6, \qquad \deg(y_1)=(-2,1), \quad \deg(y_2) = (-3, 1).
\ee
Also note that $\BP$ can be written as a GIT quotient:
\be
\BP \cong  (\C^7\setminus \{ \ud 0 \} \times \C^2\setminus \{ \ud 0 \})/(\C^\times)^2 ,
\ee
where $(\C^\times)^2 \cong \BD=\Spec (\C[\Cl(\BP)])$:
\be
(\zeta_1, \zeta_2)\cdot (x_0, \cdots, x_6, y_1, y_2) = (\zeta_1x_0, \cdots, \zeta_1x_6, \zeta_1^{-2}\zeta_2 y_1, \zeta_1^{-3}\zeta_2 y_2), \quad (\zeta_1, \zeta_2) \in (\C^\times)^2.
\ee
Then the anti-canonical divisor class $\beta=-K_{\BP}$ of $\BP$ is given by $\sum_{i=1}^9 \deg z_i = (2,2) \in \Z^2 \cong\Cl(\BP)$. Since $\BP$ is Fano, $\beta$ is ample.
Let 
\be
f(\ud z)= f(\ud x, \ud y)= y_1^2 \cdot u(\ud x) + y_2^2 \cdot v(\ud x) \in \C[z_1, \cdots, z_9],
\ee
where $u(\ud x)$ (respectively, $v(\ud x)$) is a homogeneous polynomial of degree 6 (respectively, 8) in $\C[x_0, \cdots, x_6]=\C[z_1, \cdots, z_7]$.
Then $\deg f = (2,2)=\beta$ under the identification $\Z^2 \cong \Cl(\BP)$.

Moreover, we assume that $X_{u,v}:=\{ \ud x \in \BP^6 : u(\ud x)=v(\ud x)=0 \}$ defines a smooth projective complete intersection variety of ample hypersurfaces in $\BP^6$.
Then $X_f \subset \BP$ defines a Calabi--Yau (i.e. $\deg f = \beta$) smooth ample hypersurface in $\BP$.
Note that the critical locus
\be
\operatorname{Crit}(f) =\left\{ \ud a \in \C^9 \ \middle | \    \frac{\partial  f(\ud a)}{\partial z_i} = 0, \quad i =1, \cdots, 9 \right\}
\ee
is given by
\be
\operatorname{Crit}(f) = (\C^7 \times \{ \ud 0 \}) \sqcup (\{ \ud 0 \} \times (\C^2\setminus \{ \ud 0 \}),
\ee
which is non-isolated and non-compact in $\C^9$. On the other hand, \eqref{extraisom} holds, since $X$ is even-dimensional. Hence our main theorem applies. % As we write $R(f) = \C[z_1, \cdots, z_9]/J(f)$, the Batyrev--Cox theorem yields a $\C$-vector space isomorphism:
%\be
%\bigoplus_{p=0}^{\infty} R(f)_{p \beta} =\bigoplus_{p=0}^{6} R(f)_{p \beta} \cong H_{\text{pr}}^{6}(X).
%\ee
\end{example}

\end{document}